\documentstyle[12pt,amssymb]{amsart}

\makeatletter
\def\eqalign#1{\null\,\vcenter{\openup\jot\m@th
  \ialign{\strut\hfil$\displaystyle{##}$&$\displaystyle{{}##}$\hfil
      \crcr#1\crcr}}\,}
\makeatother

      \def\g{\gamma}         
    \def\s{\sigma}  

                 \def\HD{{\mathrm{HD}}}  
         
 \def\HDess {\HD_{{\mathrm{ess}}}}

  \def\a{\alpha}        
             
\def\e{\varepsilon}      
\def\f{\varphi}            \def\d{\delta}
\def\la{\lambda}   \def\th{\vartheta}

\def\C{{\Bbb C}}         \def\Z{{\Bbb Z}}
\def\R{{\Bbb R}}             
\def\Crit{{\mathrm {Crit}}}      
      
\def\L{{\cal L}}   \def\1{{\boldsymbol 1}}

\def\and{{\mathrm and}}   \def\ov{\overline}   \def\for{{\mathrm for}}
 \def\dist{{\mathrm {dist}}}     
\def\h{{\mathrm h}}                   
\def\Vol{{\mathrm {Vol}}}        \def\diam{{\mathrm {diam}}}
\def\Const{{\mathrm {Const}}}   
\def\Comp{{\mathrm {Comp}}}   
   
    \def\Jac{{\mathrm Jac}}

\begin{document}

\title[Julia sets of holomorphic Collet--Eckmann maps]
{On measure and Hausdorff dimension of Julia sets for 
 holomorphic Collet--Eckmann maps} 
\author{Feliks Przytycki}
\address{Institute of Mathematics\\
Polish Academy of Sciences\\
ul. \'Sniadeckich 8\\
00 950 Warszawa, Poland}
\email{feliksp@@impan.impan.gov.pl}

\thanks{The author acknowledges support by Polish KBN Grant 2 P301
01307 ``Iteracje i Fraktale, II''. He expresses also his gratitude to
the MSRI at Berkeley (partial support by NSF grant DMS-9022140) and
ICTP at Trieste, where parts of this paper were written.}
\date{(version of July 26, 1995)}

\maketitle

\begin{abstract}  
Let $f:\ov\C\to\ov\C$
be a rational map 
on the Riemann sphere , such that
for every $f$-critical 
point $c\in J$ which forward trajectory does not contain 
any other critical point, 
  $|(f^n)'(f(c))|$ 
  grows exponentially fast  (Collet--Eckmann condition), 
   there are no parabolic periodic points, and else such that Julia 
   set is not the whole sphere.  Then  
smooth (Riemann) measure of the Julia set is 0.

For $f$ satisfying additionally Masato Tsujii's condition that the average 
distance of $f^n(c)$ from the set of critical points is not too small, we prove that 
Hausdorff dimension of Julia set 
is less than 2.        
This is the case for $f(z)=z^2+c$ with $c$ real, $0\in J$, 
for a positive line 
measure set of parameters $c$.
\end{abstract}

\section*{Introduction}

It is well-known that if $f:\ov\C\to\ov\C$ a rational map of the
Riemann sphere is hyperbolic, i.e. expanding on its Julia set $J=J(f)$
namely $|(f^n)'|>1$ for an integer $n>0$, then Hausdorff dimension
$\HD(J)<2$.

The same holds for a more general class of subexpanding maps, namely
such maps that all critical points in $J(f)$ are non-reccurrent,
supposed $J(f)\not=\ov\C$, see [U].

On the other hand there is an abundance of rational maps with
$J\not=\ov\C$ and $\HD(J)=2$, [Shi].

Recently Chris Bishop and Peter Jones proved that for every finitely
generated not geometrically finite Kleinian groups for the Poincar\'e
limit set $\Lambda$ one has $\HD(\Lambda)=2$. As geometrically finite
exhibits some analogy to subexpanding in the Kleinian Groups --
Rational Maps dictionary, the question arised, expressed by Ch.
Bishop and M. Lyubich at MSRI Berkeley conference in January 1995,
isn't it true for every non-subexpanding rational map with connected
Julia set, that $HD(J)=2$ ?

Here we give a negative answer. For a large class of ``non-uniformly"
hyperbolic so called Collet--Eckmann maps, studied in [P1], satisfying
an additional Tsujii condition, $ \HD(J)<2$.

\medskip\noindent {\bf Notation.} For a rational map $f:\ov\C\to\ov\C$
denote by $\Crit(f)$ the set of all critical points of $f$, i.e.points
where $f'=0$.  Let $\nu:=\sup\{\hbox{multiplicity of}f^n \hbox{at} c:
c\in\Crit(f)\cap J\}$.  Finally denote by $\Crit'(f)$ the set of all
critical points of $f$ in $J(f)$ which forward trajectories do not
contain other critical points.

We prove in this paper the following results:

\medskip\noindent {\bf Theorem A.} Let $f$ be a rational map on the
Riemann sphere $f:\ov\C\to\ov\C$, and there exist $\lambda>1, C>0$
such that for every $f$-critical point $c\in \Crit'(f)$
$$
|(f^n)'(f(c))|\ge C\lambda^n, \eqno (0.1)
$$ 
there are no parabolic periodic points, and $J(f)\not=\ov\C$.  Then
$\Vol(J(f)=0$, where Vol denotes Riemann measure on $\ov\C$.

\medskip\noindent {\bf Theorem B.} In the conditions of Theorem A
assume additionally that
$$
\lim_{t\to\infty}\limsup_{n\to\infty}{1\over n}\sum_{j=1}^n
\max(0,-\log\bigl(\dist(f^j(c),\Crit(f))\bigr)-t) = 0.  \eqno (0.2)
$$
Then $\HD(J(f))<2$.

\medskip\noindent For $f(z)=z^2+c, c\in [-2,0]$ real, it is proved in
[T] that (0.1) and (0.2) are satisfied for a positive measure set of
parameters $c$ for which there is no sink in the interval $[c,c^2+c]$.
Tsujii's condition in [T] called there {\it weak regularity} is in
fact apparently stronger than (0.2).  The set of subexpanding maps
satisfying (0.1) and weak regularity has measure 0, [T]. Thus Theorem
B answers Bishop--Lyubich's question.

\medskip\noindent {\bf Remark.} In [DPU] it is proved that for every
rational map $f\ov\C\to\ov\C, c\in\Crit'$
$$
\limsup_{n\to\infty} {1\over
n}\sum_{j=1}^n-\log\dist(f^j(c),\Crit(f))\le C_f
$$ 
where $C_f$ depends only on $f$.  Here in the condition (0.2) it is
sufficient, for Theorem B to hold, to have a positive constant instead
of 0 on the right hand side, unfortunately apparently much smaller
than $C_f$.

Crucial in proving Theorems A and B is the following intermediate
result:

\medskip\noindent {\bf Theorem 0.1 (on the existence of pacim)}, see
[P1].  Let $f:\ov\C\to\ov\C$ satisfies the assumptions of Theorem A.
Let $\mu$ be an $\a$-conformal measure on the Julia set $J=J(f)$ for
an arbitrary $\a>0$.  Assume there exists $0<\lambda'<\lambda$ such
that for every $n\ge 1$ and every $c\in \Crit'(f)$
$$
\int{d\mu \over \dist (x,f^n(c))^{(1-1/\nu)\a} }< C^{-1}(\lambda')^{\a
n/\nu}.  \eqno (0.3)
$$
Then there exists an $f$-invariant probability measure $m$ on $J$
absolutely continuous with respect to $\mu$ (pacim).

\medskip\noindent Recall that a probability measure $\mu$ on $J$ is
called $\a$-conformal if for every Borel $B\subset J$ on which $f$ is
injective $\mu(f(B)=\int_B |f'|^\a d\mu$.  In particular $|f'|^\a$ is
Jacobian for $f$ and $\mu$. The number $\a$ is called the exponent of
the conformal measure.

If $\Vol(J)>0$ then the restriction of Vol to $J$, normalized, is
2-conformal and obviously satisfies (0.3). If $\HD(J)=2$ then by [P1]
we know there exists a 2-conformal measure $\mu$ on $J$ but we do not
know whether it is not too singular, namely whether it satisfies
(0.3). Under the additional assumption (0.2) we shall prove that it is
so for every $\a$-conformal measure.

\medskip\noindent {\bf Notation.} Const will denote various positive
constants which may change from one formula to another, even in one
string of estimates.

\section{More on pacim. Proof of Theorem A.}

{\bf Proposition 1.1.} In the situation of Theorem 0.1 there exists
$K>0$ such that $\mu$-a.e.  ${dm \over d\mu} \ge K$.

\medskip\noindent {\bf Proof.} In Proof of Theorem 0.1 [P1] one
obtains $m$ as a weak* limit of a subsequence of the sequence of
measures ${1\over n}\sum_{j=0}^{n-1}f_*^j(\mu)$.

It is sufficient to prove that there exists $K>0$ and $n_0>0$ such
that for $\mu$-a.e. $y\in J(f)$
$$
{df^n_*(\mu) \over d\mu}(y)={\cal L}^n(\1) \ge K.  \eqno (1.1)
$$
Here $\cal L$ denotes the transfer operator, which can be defined for
example by ${\cal L}(\f)(y)=\sum_{f(z)=y}|f'(z)|^{-\a}\f(z)$.  $\1$ is
the constant function of value 1.  We can assume $y\notin\bigcup_{n>0}
f^n(\Crit(f))$ because
$$
\mu(\bigcup_{n>0} f^n(\Crit(f)))=0 \eqno (1.2)
$$
If a critical value for $f^n$ were an atom then a critical point would
have $\mu$ measure equal to $\infty$.

It is sufficient to prove (1.1) for $y\in B(x,\d)\cap J(f)$ for an
{\it a priori} chosen $x$ and an arbitrarily small $\d$ and next to
use the fact that there exists $m\ge 0$ such that $f^m(B(x,\d))\supset
J(f)$ (called {\it topological exactness}). Indeed
$$
\L^n(\1)(w)=\sum_{f^m(y)=w}\L^{n-m}(\1)(y)|(f^m)'|^{-\a}\ge\
(\sup|(f^m)'|)^{-\a}\L^{n-m}(y_0)
$$
where $y_0\in f^{-m}(\{w\})\cap B(x,\d)$.

Recall the estimate from [P1]. For an arbitrary $\g>1$ there exists
$C>0$ such that for every $x\in J(f)$
$$
{\cal L}^n(\1)(x)\le C+C \sum_{c\in\Crit(f)\cap J} \sum_{j=0}^{\infty}
{\g^j\la^{-j\a/\nu}\over \dist (x,f^j(f(c)))^{(1-1/\nu)\a} }.  \eqno
(1.3)
$$
By the assumptions (0.1) and (0.3) the above function is
$\mu$-integrable if $\g$ is small enough.

Pay attention to the assumption (0.3). It concerns only $c\in\Crit'$.
Fortunately there is only a finite number of summands in (1.3) for
which $f^{j_0}(c)\in\Crit, j_0\ge j$ . Each summand is integrable
because up to a constant it is bounded by $\L^j(\1)$.

So
$$\sum_{c\in\Crit(f)\cap J}
\sum_{j=s}^{\infty} {\g^j\la^{-j\a/\nu}\over \dist
(x,f^j(f(c)))^{(1-1/\nu)\a} }\to 0 \ \ \mu-{\rm a.e.} \ \hbox{as}
s\to\infty .  \eqno (1.4)
$$

Fix from now on an arbitrary $x\in J(f)$ for which (1.4) holds,
$(dm/d\mu) (x)\ge 1$ and $x\notin\bigcup_{n>0}\f^n(\Crit(f))$
(possible by (1.2) and by $\int(dm/d\mu)d\mu=1$).

We need now to repeat from [P1] a part of Proof of Theorem 0.1:

For every $y\in B(x,\d)$ and $n>0$
$$
\eqalign{ \L^n(\1)(y)&=\sum_{y'\in f^{-n}(y), {\rm
regular}}|(f^n)'(y')|^{-\a} + \sum_{(y',s) {\rm
singular}}\L^{n-s}(\1)(y')|(f^s)'(y')|^{-\a}\cr &=\sum_{{\rm
reg},y}+\sum_{{\rm sing},y} .} \eqno (1.5)
$$

We shall recall the definitions of {\it regular} and {\it singular}:
Take an arbitrary subexponentially decreasing sequence of positive
numbers $b_j, j=1,2\dots$ with $\sum b_j=1/100$.  Denote by $B_{[k}$ the
disc $B(x,(\prod_{j=1}^k (1-b_j))2\d)$. We call $s$ the {\it
essentially critical time} for a sequence of compatible components
$W_j=\Comp f^{-j}(B_{[j})$, compatible means $f(W_j)\subset W_{j-1}$,
if there exists a critical point $c\in W_s$ such that $f^s(c)\in
B_{[s}$.

We call $y'$ {\it regular} in (1.5) if for the sequence of compatible
components $W_s, s=0,1,\dots,n, W_n\ni y'$ no $s<n$ is essentially
critical.

We call a pair $(y',s)$ {\it singular} if $f^s(y')=y$ and for the
sequence of compatible components $W_j$, for $j=0,1,\dots,s$, with
$W_s\ni y'$, the integer $s$ is the first (i.e., the only) essentially
critical time.

If $\d$ is small enough then all $s$ in $\sum_{{\rm sing},x}$ are
sufficiently large that $\sum_{{\rm sing},x}\le 1/2$. This follows
from the estimates in [P1]; here is the idea of the proof:
Transforming $\sum_{{\rm sing},x}$ in (1.5) using the induction
hypothesis (1.3) we obtain the summands
$$
C{\g^j\la^{-j\a/\nu}\over \dist (x,f^{s+j-1}(f(c)))^{(1-1/\nu)\a} },\
j=0,\dots,n-s
$$ 
multiplied by
$$
\Const |(f^{s-1})'(x')|^{-\a/\nu} a_s< \g^{s-1}\la^{-(s-1)\a/\nu}.
$$
The numbers $a_s$ are constants arising from distortion estimates,
related to $b_s$. The numbers $\g^s$ swallow them and other constants.

(There is a minor inaccuracy here: $(s,x')$ is a singular pair where
the summand appears, provided $x'$ is not in the forward trajectory of
another critical point, otherwise one moves back to it, see [P1] for
details.)

Now $\sum_{{\rm sing},x}\le 1/2$ follows from (1.4).

The result is that $\sum_{{\rm reg},x}\ge 1/2$.  So by the uniformly
bounded distortion along regular branches of $f^{-n}$ on $B(x,\d)$ we
obtain
$$
\L^n(\1)(y)\ge \sum_{{\rm reg},y}\ge\Const \sum_{{\rm
reg},x}\ge\Const>0
$$
The name {\it regular} concerned formally $y'\in f^{-n}(y)$ but in
fact it concerns the branch of $f^{-n}$ mapping $y$ to $y'$ not
depending on $y\in B(x,\d)$.

By distortion of any branch $g$ of $f^{-n}$ on a set $U$ we mean
$$\sup_{z\in B}|g'(z)|/\inf_{z\in B}|g'(z)|.$$ This proves Proposition
1.1.\hfill$\square$

\medskip\noindent {\bf Corollary 1.2} In the situation of Theorem 0.1
for measure-theoretic entropy $\h_m(f)>0$.

\medskip\noindent {\bf Proof.} Denote $dm/d\mu$ by $u$.

Consider an open set $U\subset \ov\C$ intersecting $J(f)$ such that
there exist two branches $g_1$ and $g_2$ of $f^{-1}$ on it. Then by
the $f$-invariance of $m$ we have $\Jac_m(g_1)+\Jac_m(g_2) \le 1$
($=1$ if we considered all branches of $f^{-1}$). $\Jac_m(g_i)$ means
Jacobian with respect to $m$ for $g_i$.

We have $m(U)>0$ because $\mu$ does not vanish on open sets in $J$ (by
the topological exactness of $f$ on $J$) and by Proposition 1.2.  At
$m$-a.e. $x\in U$
$$
\Jac_m(g_i)(x)=u(g_i(x))|g_i'(x)|u(x)^{-1}>0,
$$
(we used here also (1.4)).

Hence $\Jac_m(g_i)<1$, so $\Jac_m(f)>1$ on the set $g_i(U), i=1,2$ of
positive measure $m$.  Now we use Rokchlin's formula and obtain
$$
\h_m(f)=\int\log \Jac_m(f) dm >0 .\eqno\square
$$

Let $\chi_m=\int\log|f'|dm$ denote characteristic Lyapunov exponent.

\medskip\noindent {\bf Corollary 1.3} In the situation of Theorem 0.1,
$\chi_m>0$.

\medskip\noindent {\bf Proof.} This Corollary follows from Ruelle's
inequality $\h_m(f)<2\chi_m$, see [R].

\medskip\noindent {\bf Proof of Theorem A.} Suppose $\Vol(J(f))>0$.
After normalization we obtain a 2-conformal measure $\mu$ on $J(f)$
and by Theorem 0.1 and Corollary 1.3 a pacim $m$ with $\chi_m>0$.  By
Pesin's Theory [Pesin] in the iteration in dimension 1 case [Le] ([Le]
is on the real case, but the complex one is similar), for $m$-a.e.
$x$, there exists a sequence of integers $n_j\to\infty$ and $r>0$ such
that for every $j$ there exists a univalent branch $g_j$ of $f^{-n_j}$
on $B_j:=B(f^{n_j}(x),r)$ mapping $f^{n_j}(x)$ to $x$ and $g_j$ has
distortion bounded by a constant.  By $\chi_m>0 \ \diam
g_j(B(f^{n_j}(x),r)\to 0$. (This follows also automatically from the
previous assertions by the definition of Julia set [GPS].) Now we can
forget about the invariant measure $m$ and go back to $\Vol$.  Because
$J(f)$ is nowhere dense in $\ov\C$, there exists $\e>0$ such that for
every $z\in J(f)$
$$
{\Vol (B(z,r)\setminus J(f)) \over \Vol (B(z,r))}>\e.
$$

Bounded distortion for $g_j$ on $B(z,r)$, $z=f^{n_j}(x)$ allows to
deduce that the same part of each small disc$\approx g_j(B_j)$ around
$x$ is outside $J(f)$ , up to multiplication by a constant.  This is
so because we can write for every $X\subset B(z,r), y\in B(z,r)$
$$
\Vol g_j(X)\approx |g_j'(y)|^2 \Vol(X) \eqno (1.6)
$$ 
where $\approx$ means up to the multiplication by a uniformly bounded
factor.  So $x$ is not a density point of $J(f)$.  On the other hand
a.e. point is a density point. So $\Vol J(f)=0$ and we arrived at a
contradiction.\hfill$\square$

\section{Proof of Theorem B.}

As mentioned in the Introduction, the following result is crucial:

\medskip\noindent {\bf Lemma 2.1.} Under the conditions of Theorem B,
i.e. in the situation of Theorem A and assuming (0.2), the condition
(0.3) holds for every $\a$-conformal measure, $\a>0$.

\medskip\noindent {\bf Proof.\ Step 1.} Denote the expression from
(0.2)
$$
\max\bigl(0,-\log\inf_{c\in\Crit'(f)}\dist(f^n(c),\Crit(f))-t\bigr)
$$
by $\f_t(n)$.  Consider the following union of open-closed intervals
$$
A_t':=\bigcup_n (n,n+\f_t(n)\cdot K_f] \ \hbox{and write} \
A_t:=\Z_+\setminus A_t',
$$
for an arbitrary constant $K_f> \nu/\log\la$.

By (0.2) for every $a>0$ there exist $t>0$ and $n(a,t)$ such that for
every $n\ge n(a,t)$
$$
A_t\cap [n,n(1+a)]\not=\emptyset \eqno (2.1)
$$

Moreover, fixed an arbitrary integer $M>0$, we can guarantee for every
$n'\ge n(1+a),n\ge n(a,t)$
$$
\sharp (A_t \cap \{j\in [n,n']: j \hbox{divisable by}M\}) \ge {1\over
2M} (n'-n).  \eqno (2.2)
$$

Observe that for every $n_0,n$, (2.1) transforms into
$$
[n_0+n,n_0+n+a(n_0+n)]=[n_0+n,n_0+n+a({n_0 \over n}+1)n].
$$
The result is that if $n \ge bn_0$ for an arbitrary $b>0$ then
$$
A_t\cap [n_0+n,n_0+n+a(b^{-1}+1)n]\not=\emptyset. \eqno (2.3)
$$
Denote in the sequel $a(b^{-1}+1)$ by $a'$.

\medskip\noindent {\bf Step 2.} Observe now that if $n\in A_t$ then
for every $c\in\Crit'(f)$ there exist branches $g_s, s=1,2,\dots,n-1$ of
$f^{-s}$ on $B_n:=B(f^n(c),\d)$ with uniformly bounded distortions,
where $\d=\e\exp-t\nu$ for a constant $\e$ small enough. Sometimes to
exhibit the dependence on $n$ we shall write $g_{s,n}$.

Indeed, define $g_s$ on $B_{[s}=B(f^n(c),\prod_{j=1}^s (1-b_j)2\d)$
for $s=1,2,\dots,n-1$ according to the procedure described in Proof of
Theorem A.  If there is an obstruction, namely $s$ an essential
critical time, then for every $z\in B_{[s}$
$$
|g_{s-1}'(z)| \le \la^{-s}\th^s \le \exp (-s\nu/K_f) \eqno (2.4)
$$
for $\th>1$ arbitrarily close to 1 (in particular such that
$K_f>{\nu\over \log\la-\log\th}$) and for $s$ large enough. The
constant $\th$ takes care of distortion. (2.4) holds for $z=f^s(q)$,
where $q$ is the critical point making $s$ critical time, without
$\th$ by (0.1) (with the constant $C$ instead). The small number $\e$
takes care of $s$ small, which cannot be then essential critical.

The inequality (2.4) and rooting ($1/\nu$ to pass from $s-1$ to $s$)
imply $\f_t(f^{n-s}(c))\ge s/K_f$, so $n\notin A_t$, a contradiction.

\medskip\noindent {\bf Step 3.} By uniformly bounded distortion for
the maps $g_{j,n}$, $n\in A_t$ we obtain (compare (1.6)) for every
$n_0>0$ large enough, $ c\in\Crit'$
$$
\mu B(f^{n_0}(c), r_j)\approx r_j^\a \eqno (2.5)
$$
for a sequence $r_j, j=1,2,\dots$ such that
$$
r_1>\exp -Lbn_0 \eqno (2.6)
$$
$$
r_{j+1} > r_{j}^{1+ \s} \eqno (2.7)
$$
$$
\and \ r_{j+1}<r_{j}/2.  \eqno (2.8)
$$
Here $L:=2\sup |f'|$ and $b, \s$ are arbitrarily close to 0.

Indeed, we can find $r_j$ satisfying the conditions above by taking
$$
r_j:=\diam g_{n_j, n_0+n_j}(B(f^{n_0+n_j}(c),\d))
$$ 
where $n_j\in A_t$ are taken consecutively so that
$$
n_{j+1}\in [n_j+(1+\th)n_j, n_j+(1+\th)(1+a')n_j \ \for \ j\ge 2 \
\and
$$
$$
n_1\in [n_0+bn_0, n_0+bn_0+a'b n_0],
$$
where $\th>0$ is an arbitrary constant close enough to 0.

This gives
$$
r_{j+1}/r_j \ge \exp (-2(\log L)a'n_j) \eqno (2.9)
$$
To conclude we need to know that $r_j$ shrink exponentially fast with
$n_j\to\infty$. For that we need the following fact (see for example
[GPS], find the analogous fact in Proof of Theorem A):

\medskip\noindent {\bf (*)} For every $r>0$ small enough and $\xi,
C>0$ there exists $m_0$ such that for every $m\ge m_0, x\in J(f)$ and
branch $g$ of $f^{-m}$ on $B(x,r)$ having distortion less than $C$ we
have $\diam g(B(x,r))<\xi r$.

\medskip\noindent Apply now (2.2) to $n=n_0, n'=n_j+n_0$. We obtain a
``telescope": For all consecutive $\tau_1,\tau_2,\dots,\tau_{k(j)}\in
A_t\cap [n_0, n_j+n_0]$ divisible by $M$
$$
g_{\tau_{i+1}-\tau_i, \tau_{i+1}} (B(f^{\tau_{i+1}}(c),\d))\subset
B(f^{\tau_i}(c),\d/2)
$$
for $M\ge m_0$ from (*).

Hence using (2.2)
$$
r_j\le 2^{-n_j/2M} .  \eqno (2.10)
$$
Denote $2a'\log L$ by $\g$ and $(\log 2)/2M$ by $\g'$. (2.9) and
(2.10) give
$$
r_{j+1} \ge r_j\exp -\g n_j \ge r_j (\exp -\g'n_j)^{\g/\g'} \ge
r_j^{1+\g/\g'}.
$$
As $\g'$ is a constant and $\g$ can be done arbitrarily small if $a$
is small enough, we obtain (2.7).

The condition (2.8) follows from the fact that for $n_0$ large enough
all $n_{j+1}-n_j$ are large enough to apply (*).

\medskip\noindent {\bf Conclusion.} We obtain
$$
\int{d\mu\over\dist(x,f^{n_0}(c))^{(1-1/\nu)\a}} \le
$$
$$
\mu(\ov\C \setminus B(f^{n_0}(c), r_1)) {1\over r_1^{(1-1/\nu)\a}} +
\sum_{j\ge 2} \mu(B(f^{n_0}(c), r_{j-1})\setminus B(f^{n_0}(c), r_j))
{1\over r_j^{(1-1/\nu)\a}} \le
$$
$$
 \exp (Lbn_0(1-1/\nu)\a) +\Const\sum_{j\ge 2} {r_{j-1}^\a \over
r_j^{(1-1/\nu)\a}} \le
$$
$$
(\exp (Lb(1-1/\nu)\a))^{n_0} + \Const\sum_{j\ge 2} r_{j-1}^\a
r_{j-1}^{-(1-1/\nu)\a (1+\s)}.
$$

The latter series has summands decreasing exponentially fast for $\s$
small enough so it sums up to a constant, hence the first summand
dominates. We obtain the bound by $(\la')^{n_0}$ with $\la'>1$
arbitrarily close to 1. Thus (0.3) has been proved.\hfill$\square$

\medskip\noindent {\bf Remark 2.2.} The only result in our disposal on
the abundance of non-subexpanding maps satisfying (0.1) and (0.2) is
Tsujii's one concerning $z^2+c, c$ real (see Introduction). For this
class however the exponential convergence of 
$$\diam\Comp
f^{-n_j}(B(f^{n_j+n_0}(0),\d)$$ to 0 follows from [N] (the component
containing $f^{n_0}(c)$).  So restricting our interests to this class
we could skip (2.2) and considerations leading to (2.10) above.

By [N] $\diam\bigl(\Comp( f^{-n}(B(x,\d)))\cap \R\bigr) <
C\tilde\la^{-n}$ for some constants $C>0, \tilde\la>1, \d$ small
enough and every component Comp. Just the uniform convergence of the
diameters to 0 as $n\to\infty$ follows from [P1], but I do not know
how fast is it.

\medskip\noindent {\bf Proof of Theorem B.} Suppose that $\HD(J)=2$.
Then there exist a 2-conformal measure $\mu$ on $J$. This follows from
the existence of an $\a$-conformal measure for $\a=\HDess(J)$, where
$\HDess$ is the essential Hausdorff dimension which can be defined for
example as supremum of Hausdorff dimensions of expanding Cantor sets
in $J$, see [DU] [P2] and [PUbook], and from $\HDess(J)=\HD(J)$, see
[P1]. The former holds for every rational map, the latter was proved
in [P1] only for Collet--Eckmann maps.

By Lemma 2.1 the condition (0.3) holds, hence there exists a pacim
$m\ll\mu$.  Moreover $\chi_m>0$ by Corollary 1.3. As in Proof of
Theorem A by Pesin Theory there exists $X\subset J$, $m(X)=\mu(X)=1$,
such that for every $x\in X$ there exists a sequence of integers
$n_j(x)\to\infty$, $r>0$ and univalent branches $g_j$ of $f^{-n_j}$ on
$B(f^{n_j}(x),r)$ mapping $f^{n_j}$ to $x$ with uniformly bounded
distortion. Write $B_{x,j}:=g_j (B(f^{n_j}(x),r))$.

We obtain for every $x\in X$ applying (1.6) to Vol and $\mu$
(similarly as in Proof of Theorem A)
$$
\mu(B_{x,j})\le\Const\Vol(B_{x,j})\le\Const(\Vol B(x,\diam B_{x,j})) .
$$
If $\Vol X=0$ then there exists a covering of $X$ by discs
$B(x_t,\diam B_{x_t,j_t}), t=1,2,\dots$ which union has $\Vol<\e$ for
$\e$ arbitrarily close to 0, of multiplicity less than a universal
constant (Besicovitch's theorem). Hence
$$
\e\ge\Const\sum_t\Vol B(x_t,\diam B_{x_t,j_t})\ge\Const \mu \sum_t
B_{x_t,j_t}\ge 1,
$$
a contradiction. Hence $\Vol J\ge \Vol X>0$.

This contradicts Theorem A that $\Vol J=0$ and the proof is over.

Remark that we could end the proof directly: As in Proof of Theorem A
we show that no point of $X$ is a point of density of the $\Vol$
measure.  Hence $\Vol X=0$. (I owe this remark to M.
Urba\'nski.)\hfill$\square$

\end{document}